\newif\ifmakeplots
\makeplotstrue
\newif\ifspringer
\springerfalse

\ifspringer
    \documentclass[smallextended,envcountsect]{svjour3} 
\else
    \documentclass[12pt,notitlepage]{article}
\fi

\newcommand{\datadir}{data}

\ifspringer
\else
    \usepackage[margin=1.25in, headheight=42pt]{geometry}
    \usepackage{amsthm}
\fi
\usepackage{fancyhdr}
\usepackage{amsmath,thmtools}
\usepackage{algorithm, algorithmic}
\usepackage{pgfplots}
\usepackage[hidelinks]{hyperref}
\usepackage[style=alphabetic,backend=bibtex,giveninits=true,url=false,doi=false,natbib=true,maxbibnames=99]{biblatex}
\def\aspectratio{1.3}

\newcommand{\ie}{\emph{i.e.}}
\newcommand{\eg}{\emph{e.g.}}
\newcommand{\reals}{\mathbf R}

\newcommand{\dom}{\mathop{\bf dom}}

\newcommand{\symm}{\mathbf S}

\newcommand{\argmin}{\mathop{\rm argmin}}
\newcommand{\card}{\mathop{\bf card}}
\newcommand{\diag}{\mathop{\bf diag}}
\newcommand{\prox}{{\bf prox}}
\newcommand{\dist}{\mathbf{dist}}
\newcommand{\ones}{\mathbf{1}}

\definecolor{darkgreen}{rgb}{0, .5, 0}
\usepgfplotslibrary{groupplots}
\pgfplotsset{scaled x ticks=false}
\pgfplotsset{scaled y ticks=false}

\usepgfplotslibrary{groupplots}
\pgfplotsset{
    layers/my layer set/.define layer set={
        background,
        background1,
        background2,
        background3,
        background4,
        main,
        foreground
    }{
    },
    set layers=my layer set,
}

\ifspringer
\else
    \fancypagestyle{plain}{%
        \fancyhead{}
        \chead{\textsc{For professional clients and qualified investors only \\
               Not for public distribution \\
               This version posted with permission}}
    }
\fi

\title{Portfolio Construction as\\Linearly Constrained Separable Optimization}

\author{Nicholas Moehle \and Jack Gindi \and Stephen Boyd \and Mykel J. Kochenderfer}

\begin{document}

\maketitle

\begin{abstract}
Mean--variance portfolio optimization problems often involve separable nonconvex terms,
including penalties on capital gains, integer share constraints,
and minimum nonzero position and trade sizes.
We propose a heuristic algorithm for such problems
based on the alternating direction method of multipliers (ADMM).
This method allows for solve times in tens to hundreds of milliseconds
with around 1000 securities and 100 risk factors.
We also obtain a bound on the achievable performance.
Our heuristic and bound are both derived from similar results
for other optimization problems with a separable objective and affine equality constraints.
We discuss a concrete implementation in the case where the separable terms
in the objective are piecewise quadratic,
and we empirically demonstrate its effectiveness for tax-aware portfolio construction.
\end{abstract}

\section{Introduction}
\label{s-intro}
The mean--variance portfolio optimization problem of \citet{markowitz1952}
has a quadratic objective and linear equality constraints,
allowing for a simple analytical solution.
The problem can be extended by including position limits or a long-only constraint
\citep{markowitz1955optimization, sharpe1963simplified, grinold2000active}.
Although the resulting problem no longer has an analytical solution,
it can be efficiently solved as a quadratic program \citep[pp.~55--156]{cvxbook}.
Many extensions have been proposed that take into account investment restrictions,
risk or leverage limits, and trading costs \citep{boyd2017multi}.
When these extensions preserve the convexity of the portfolio optimization problem,
the extended problems can be solved quickly and reliably \citep[\S 1.3.1]{cvxbook}.

Some important portfolio optimization problems involve
nonconvex constraints and objective terms.
For example, we may require that asset holdings be in integral numbers of shares, or
we may want to penalize realized capital gains \citep{moehle2020tax}
or the number of securities traded \citep{lobo2007portfolio}.
These extensions are often irrelevant or negligible for large institutional accounts,
but they can be important for small accounts
such as those arising in separately managed account (SMA) platforms \citep{benidis2018optimization}.
Some of these nonconvex portfolio optimization problems can be reformulated as
mixed-integer convex optimization problems \citep{mansini2015linear} and solved exactly using commercial or open source solvers.
Such solvers are often fast, but occasionally have very long solve times,
often hundreds of times more than for similar convex problems.

This article focuses on nonconvex portfolio optimization problems
with nonconvex penalties and constraints that separate across assets.
This problem is a special case of the \emph{separable--affine problem},
\ie, the problem of minimizing a separable objective function with affine equality constraints.
We propose a heuristic algorithm for solving such problems
based on the alternating direction method of multipliers (ADMM).
This method solves problems to moderate accuracy quickly,
even when the separable functions are very complicated.
% (In portfolio optimization, this occurs when optimizing over tax liability function
% when there are many tax lots, especially when trading within the 30-day wash sale window.)
This speed is obtained because each separable function is interfaced only
through its proximal operator, which only involves a few arithmetic operations.
This separable--affine form also makes it easy to compute a lower bound on the problem value
%using the convex relaxation obtained 
by replacing the separable functions by their convex envelopes
and solving the resulting convex optimization problem.
We give a fast algorithm for constructing these convex envelopes for piecewise quadratic functions,
which often appear in portfolio optimization.

This article is structured as follows. 
In section~\ref{s-port-opt}, we introduce the portfolio optimization problem
with separable nonconvex holding and trading costs.
In section~\ref{s-sap}, we introduce the separable--affine problem
and show how the portfolio optimization problem can be expressed in this form.
In section~\ref{s-bounds}, we discuss lower bounds on the optimal value 
of the separable--affine problem.
In section~\ref{s-admm} we show how to solve separable--affine problems using ADMM.
In section~\ref{s-implementation}, we discuss implementation details 
when the separable functions are piecewise quadratic.
Section~\ref{s-example} concludes with an application to tax-aware portfolio 
construction.

\section{Portfolio optimization}
\label{s-port-opt}

% \paragraph{Cash.}
% %The initial cash held in the portfolio is $c^{\rm init}\in\reals$.
% The fraction of the portfolio held in cash is $c = 1 - \ones^T h$.
% We aim to maintain the cash fraction within some specified bounds,
% \ie, $c^{\rm lb} \le c \le c^{\rm ub}$,
% which can be expressed in terms of $h$ as
% \[
% c^{\rm lb} \le 1 - \ones^T h \le c^{\rm ub}.
% \]

\subsection{Basic portfolio optimization problem}
We first present mean--variance portfolio optimization
with separable holding and trading costs.
The problem is to decide how much to invest in each of $l$ assets.
This decision is represented by $h\in\reals^l$,
where $h_i$ is the fraction of the portfolio value to be invested in asset $i$.

We choose $h$ by solving the optimization problem
\begin{equation}\label{e-port-opt}
\begin{array}{ll}
\text{maximize} & \alpha^T h - \gamma^{\rm risk} h^T V h
                    - \gamma^{\rm trd} \phi^{\rm trd}(h - h^{\rm init}) - \gamma^{\rm hld}\phi^{\rm hld}(h) \\
\text{subject to} & \eta^{\rm lb} \le \ones^T h \le \eta^{\rm ub}.
\end{array}
\end{equation}
The objective of \eqref{e-port-opt} trades off expected return, 
risk, trading costs, and holding costs,
with positive tradeoff parameters 
$\gamma^{\rm risk}$, $\gamma^{\rm trd}$, and $\gamma^{\rm hld}$.
The vector $\alpha\in\reals^l$ is the expected return forecast for the $n$ assets,
meaning $\alpha^T h$ is the expected portfolio return.
The matrix $V\in\symm_{++}^l$
(the set of symmetric positive definite $l \times l$ matrices)
is the asset return covariance matrix,
making $h^T V h$ the variance of the portfolio return.
We assume $V$ has the traditional factor model form
\begin{equation}\label{e-factor-model}
V = X\Sigma X^T + D,
\end{equation}
where $X\in \reals^{n \times k}$ is the factor exposure matrix,
$\Sigma\in \symm_{++}^k$ is the factor covariance
matrix, and $D\in\symm^{n}_{++}$ is the diagonal matrix of idiosyncratic 
variances with $D_{ii} > 0$ \citep{grinold2000active, boyd2017multi}.

The constraint specifies that the fraction of the account value that is invested
(given by $\ones^T h$)
is between $\eta^{\rm lb}$ and $\eta^{\rm ub}$.
Setting $\eta^{\rm lb} = 0.98$ and $\eta^{\rm ub} = 0.99$, for example,
means that between 98\% and 99\% of the account value must be invested
after the trade is carried out, leaving 1--2\% in cash.

The trading cost $\phi^{\rm trd}(h - h^{\rm init})$ is
the cost of trading from the initial portfolio $h^{\rm init}$ to $h$,
and $\phi^{\rm hld}(h)$ is the cost of holding portfolio $h$.
We do \emph{not} assume these two functions are convex,
but we do assume they are separable,
\ie,
\[
\phi^{\rm trd}(u) = \sum_{i=1}^l \phi^{\rm trd}_i(u_i),
\qquad
\phi^{\rm hld}(h) = \sum_{i=1}^l \phi^{\rm hld}_i(h_i).
\]

Because $\phi^{\rm trd}$ and $\phi^{\rm hld}$ are not convex,
problem~\eqref{e-port-opt} is difficult to solve exactly in general.
One way to solve it is to reformulate
it as a mixed-integer convex problem,
typically a mixed-integer quadratic program or a mixed-integer second-order cone program,
both of which can be solved using standard methods.
In practice, these methods often solve problem instances quickly,
but sometimes the solution times can be extremely long.

In this paper, we propose a different approach,
which exploits the fact that the nonconvex terms are separable.
Problems with separable nonconvex terms have been studied extensively,
and many effective solution methods have been proposed for them.
We discuss problems of this type in more detail in section~\ref{s-sap}.

\subsection{Examples of separable trading costs}
We give some examples of trading cost functions $\phi^{\rm trd}$.
Although we discuss each individually,
in practice these functions would be combined together
into a single composite trading cost function.

\paragraph{Transaction cost.}
The traditional transaction cost model is
\[
\phi^{\rm trd}_i(u) = s_i |u| + d_i |u|^{3/2}.
\]
The first term models the cost of crossing the bid--ask spread,
where $s_i \ge 0$ is one-half the bid-ask spread of asset $i$,
and $|\cdot|$ denotes the element-wise absolute value.
The second term models the cost of price impact,
where $d_i \ge 0$ is the market impact parameter for asset $i$.
The $3/2$ power is applied elementwise.

\paragraph{Minimum trade size.}
A minimum trade size is a constraint that for each asset $i$,
we either do not trade it, or we trade at least $u_i^{\rm min}$ of it.
In this case, we have
\[
\phi^{\rm trd}_i(u_i) = \begin{cases}
0 & \text{if $|u_i| \ge u_i^{\rm min}$ or $u_i = 0$}  \\
\infty & \text{otherwise.} \\
\end{cases}
\]

\paragraph{Per-trade cost.}
A per-trade cost has the form
\[
\phi^{\rm trd}_i(u_i) = \begin{cases}
0 & \text{if $u_i = 0$} \\
c^{\rm trd}_i & \text{otherwise,}
\end{cases}
\]
where $c^{\rm trd}_i$ is the cost of trading asset $i$.
If these costs are the same for all assets,
the portfolio-level trading cost is proportional to the cardinality $\card(u)$,
\ie, number of nonzero elements in $u$.
This can be used to model a per-trade commission levied on all trades.

\paragraph{Tax liability.}
Let $L: \reals^l \to \reals$ denote the tax liability function,
such that $L(h - h^{\rm init})$ is the immediate tax liability incurred from realized capital gains
by trading into portfolio $h$.
This function is separable across the assets, \ie, it has the form
\[
\phi^{\rm trd}_i(u_i) = L_i(u_i),
\]
where $L_i(u_i)$ is the tax liability from trading asset $i$.
It is also piecewise affine and has domain $[-h^{\rm init}_i, \infty)$.
An explicit description of $L_i$ is complicated \citep[\S 3]{moehle2020tax}.
(The function $L_i$ is especially complicated when
the asset $i$ has been traded within the past $30$ days because of the 
wash sale rule, but it is always a piecewise affine function.)

\subsection{Examples of separable holding costs}
In this section, we provide some examples of trading cost 
functions $\phi^{\rm hld}$.
Although we discuss each individually, they are combined
into a composite holding cost function.

\paragraph{Position limits.}
Asset-level position limits have the form
\begin{equation}
\label{e-position-limits}
    \phi^{\rm hld}_i(h_i) = \begin{cases}
    0 & \text{if $h^{\rm lb}_i \le h_i \le h^{\rm ub}_i$} \\
\infty & \text{otherwise,}
\end{cases}
\end{equation}
where limits for asset $i$ satisfy $h^{\rm lb}_i \le h^{\rm ub}_i$.

\paragraph{Minimum holding size.}
A minimum holding size is a constraint that for each asset $i$,
we either do not hold it, or we hold at least $h^{\rm min}_i$ of it.
In this case, we have
\[
\phi^{\rm hld}_i(h_i) = \begin{cases}
0 & \text{if $|h_i| \ge h^{\rm min}_i$ or $h_i = 0$}  \\
\infty & \text{otherwise.} \\
\end{cases}
\]

\paragraph{Per-asset holding cost.}
A per-asset holding cost has the form
\begin{equation}
\label{e-holding-cost}
\phi^{\rm hld}_i(h_i) = \begin{cases}
0 & \text{if $h_i = 0$} \\
c^{\rm hld}_i & \text{otherwise,}
\end{cases}
\end{equation}
where $c^{\rm hld}_i$ is the cost of holding asset $i$.
If these costs are the same for all assets,
the portfolio-level holding cost is proportional to the cardinality $\card(h)$,
\ie, the number of nonzero elements in $h$.
Penalizing this value models the account overhead associated
with maintaining a portfolio with many assets.

\paragraph{Integer share constraint.}
We can restrict the portfolio to hold an integer number of shares of each asset,
\ie, $h_i / p_i \in \mathbf Z$,
where $p_i$ is the per-share price of asset $i$.
In this case, we have
\begin{equation*}
\phi^{\rm hld}_i(h_i) = \begin{cases}
0 & \text{if $h_i/p_i \in\mathbf Z$} \\
\infty & \text{otherwise.} \\
\end{cases}
\end{equation*}

\section{Separable--affine problem}
\label{s-sap}
We propose to solve the portfolio optimization problem~\eqref{e-port-opt}
as a \emph{separable--affine problem} (SAP).
In this section, we introduce the SAP problem and discusses some of its properties.

\paragraph{Definition.}
The separable--affine problem is
\begin{equation}
  \label{e-sap}
  \begin{array}{ll}
    \mbox{minimize} & \sum_{i=1}^n f_i(x_i)\\
    \mbox{subject to} & Ax = b,
  \end{array}
\end{equation}
with variable $x\in \reals^n$.
The parameters are $A \in \reals^{m\times n}$ and $b \in \reals^m$,
as well as the separable functions $f_i:\reals \to \reals \cup \{\infty\}$.
The (separable) objective is $f(x) = \sum_{i=1}^n f_i(x_i)$,

We use infinite values of $f_i$ to encode constraints, and define
$\dom f_i = \{ x_i \mid f_i(x_i) < \infty \}$
and $\dom f = \dom f_1 \times \cdots \times \dom f_n$.
We will assume that for each $i$, $\dom f_i$ is a non-empty union of a finite number of intervals.
We also assume each $f_i$ is closed.
We say that $x$ is feasible if $x \in \dom f$ and $Ax=b$,
and define the (equality constraint) residual associated with $x$
as $r = Ax-b$.
We denote the optimal value of the SAP \eqref{e-sap} as $p^\star$
and a solution (if one exists) as $x^\star$.

\paragraph{Scaling.} \label{s-prob-scaling}
We observe here for future use that both the variables and constraints in the
SAP \eqref{e-sap} can be scaled, yielding another (equivalent) SAP.
Let $E\in \reals^{n \times n}$ be diagonal and invertible, and
$D\in \reals^{m \times m}$ be invertible.
With the change of variable $\tilde x = E^{-1}x$, and scaling the
equality constraints by the matrix $D$, we obtain the problem
\begin{equation}
  \label{e-sep-probs-scaled}
  \begin{array}{ll}
    \mbox{minimize} & \sum_{i=1}^n f_i(E_{ii}\tilde x_i)\\
    \mbox{subject to} & DAE \tilde x = Db,
  \end{array}
\end{equation}
with variable $\tilde x$, which is also an SAP with data
\[
\tilde A = DAE, \qquad
\tilde b = Db, \qquad
\tilde f_i(\tilde x_i) = f_i(E_{ii}\tilde x_i), \quad i=1, \ldots, n.
\]
From a solution $\tilde x^\star$ of this problem, we can recover a solution
of the original SAP as $x ^\star = E\tilde x^\star$.

\subsection{Portfolio optimization as a SAP}
To express problem~\eqref{e-port-opt} as a separable--affine problem,
we introduce two additional variables:
the cash fraction $c = 1 - \ones^T h$ and the factor exposure vector $y = C^TXh$,
where $C$ is a Cholesky factor of $\Sigma$, \ie, $CC^T = \Sigma$.
Problem~\eqref{e-port-opt} is then
\begin{equation}
    \label{e-app-port-orig}
    \begin{array}{ll}
      \mbox{maximize} & \alpha^T h -  \gamma^{\rm risk} (y^T y + h^T D h)
              - \gamma^{\rm trd} \phi^{\rm trd}(h - h^{\rm init}) - \gamma^{\rm hld} \phi^{\rm hld}(h) \\
      \mbox{subject to}
          & y = C^TXh \\
          & c + \ones^T h = 1 \\
          & 1 - \eta^{\rm ub} \leq c \leq 1 - \eta^{\rm lb}. \\
    \end{array}
\end{equation}
The variables are $h\in\reals^l$, $c\in\reals$, and $y\in\reals^k$.
Minimizing the negative of the objective yields a separable--affine problem;
the exact values of the parameters $f$, $A$, and $b$ are given in appendix~\ref{s-portfolio-params}.

\subsection{Solving SAPs}
\paragraph{Convex case.}
If the functions $f_i$ are all convex,
the SAP~\eqref{e-sap} is a convex optimization problem that is readily solved.

\paragraph{Exhaustive search.} If the number of degrees of freedom $n-m$ is very small (say, no more than~3 or~4),
we can solve the SAP by exhaustive search, which involves
evaluating $f(x)$ on a grid of points in the subspace $\{x \in\reals^n \mid Ax = b\}$.

\paragraph{Divide and conquer.}
When $m$ is very small (say, no more than~3), we can use
variations on dynamic programming or divide and conquer to solve the SAP.
For a subset $S \subseteq \{1, \ldots, n\}$, we define the value function
\[
V_S(z) = \inf \bigg\{ \sum_{i \in S} f_i(x_i) \;\bigg|\; \sum_{i \in S} x_i a_i = z \bigg\},
\]
where $a_i$ are the columns of $A$.
When $m$ is very small, we can represent $V_S$ by its values on a grid, or
by the coefficients in a suitable basis.

We have
$V_{\{k\}}(z) = f_k(x_k^\star)$ if $x_k^\star a_k = z$ for some $x_k^\star$,
and $V_{\{k\}}(z) = \infty$ otherwise.
We also have $V_{\{1, \ldots, n\}}(b) = p^\star$,
as well as the dynamic programming property
\[
V_{S \cup T}(z) = \inf_{u,v} \{ V_S(u) + V_T(v) \mid u+v = z \},
\]
for disjoint $S$ and $T$.
When $m$ is very small, we can form $V_{S\cup T}$ by brute force search over
a grid of values $u$ in $\reals^m$.
Thus, we can form $V_{S \cup T}$, given $V_S$ and $V_T$
(modulo our gridding or basis approximation).

This suggests the following divide and conquer method.
We start with the collection of $n$ value functions $V_S$ with $S$ a singleton
(which are the same as the functions $f_i$).
We then combine pairs, using the combining formula above, to obtain
around $n$ value functions $V_S$, with $|S|=2$.
We continue this way around $\log_2 n$ times to obtain $V_{\{1, \ldots, n\}}$.
Evaluating this function at $z=b$ gives $p^\star$.

Examining the portfolio problem~\eqref{e-port-opt}
(or its separable--affine form \eqref{e-app-port-orig}),
note that $m = k + 1$, where $k$ is the number of factors in the risk model.
This implies that portfolio optimization problems with separable nonconvexities
can be solved exactly if the number of risk factors is small.

\paragraph{The general case is hard.}
SAP includes mixed-integer linear programs
as a special case, which in turn includes as a further special case
the 3-SAT problem, which is NP-complete.
Thus, we do not expect to develop a global optimization method that is efficient in
the worst case; we only hope for efficiency of a global method
on many (or even just some) problem instances that arise in practice.

\subsection{Special cases}

The SAP includes a number of well-known problem
classes as special cases.

\paragraph{LP and MILP.}
With $f_i(x_i) = c_i x_i + I_+(x_i)$, where $I_+$ is the indicator function
of $\reals_+$,
\[
I_+(u) = \begin{cases}
0 & \text{if $u \ge 0$} \\
\infty & \text{if $u < 0$,}
\end{cases}
\]
the SAP reduces to the generic linear programming (LP) problem.
Adding the indicator function of $\{0,1\}$ to $f_i$ we obtain a general mixed-integer
linear program (MILP).

\paragraph{Indefinite quadratic programming.}
The SAP generalizes the indefinite quadratic programming problem (IQP)
\begin{equation}
  \label{e-qp}
  \begin{array}{ll}
    \mbox{minimize} & x^T P x + q^T x \\
    \mbox{subject to} & Ax = b \\
                      & x \geq 0,
  \end{array}
\end{equation}
where $P$ is a symmetric matrix.
We factor $P$ as $P = FDF^T$, where $D$ is $r \times r$ diagonal with
nonzero diagonal entries and $F \in \reals^{n \times r}$.
This can be obtained from an eigendecomposition of $P$.
Thus, we have $x^TPx = \sum_{i=1}^n D_{ii} z_i^2$ for $z = F^T x$.
Adding this new variable, we obtain the SAP
\[
  \begin{array}{ll}
    \mbox{minimize} &
\sum_{i=1}^n (q_i (x_i) + I_+(x_i)) + \sum_{i=1}^r D_{ii} z_i^2 \\
    \mbox{subject to} & \left[ \begin{array}{cc} A & \phantom{-}0 \\ F^T & -I \end{array}
\right] \left[ \begin{array}{c} x \\ z \end{array}\right] =
\left[ \begin{array}{c} b \\ 0 \end{array}\right],
  \end{array}
\]
with variable $(x,z)$.
This is equivalent to the IQP above.
Adding the indicator function of $\{0,1\}$ to each $f_i$ yields
a general mixed-integer IQP.

\paragraph{Limitations.}
The feasible set of a SAP
consists of the Cartesian product of unions of intervals
(\ie, $\dom f_i$) and the affine set $\{x \mid Ax=b\}$.
The feasible set of an SAP is therefore the union of a finite number of convex polyhedra. 
This observations tells us that optimization problems with non-polyhedral feasible sets,
such as problems with quadratic constraints, cannot be exactly represented as SAPs.

\section{Bounds}
\label{s-bounds}
In this section, we describe a lower bound on $p^\star$
obtained by solving a convex relaxation of \eqref{e-sap}.

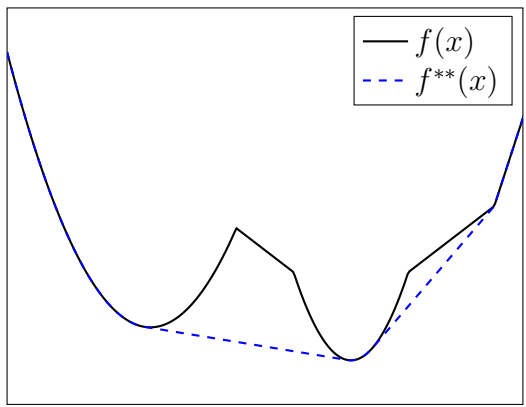
\begin{figure}
\centering
\ifmakeplots
    \begin{tikzpicture}
        \begin{axis}[xmin=-1,
                     xmax=8,
                     ymin=-4,
                     width = 0.45\textwidth,
                     height = 0.45\textwidth / \aspectratio,
                     scale only axis=true,
                     ymax=5,
                     ticks=none,
                     legend cell align={left}
                     ]

            % Nonconvex function.
            \addplot [thick, black] table{\datadir/example_function.csv};

            % Convex envelope of nonconvex function.
            \addplot [thick, dashed, blue] table{\datadir/example_function_envelope.csv};

            \legend{$f(x)$, $f^{**}(x)$}

        \end{axis}
    \end{tikzpicture}
\fi
\caption{A nonconvex function $f$ (solid black line)
and its convex envelope $f^{**}$ (dashed blue line).}
\label{f-conv-env-example}
\end{figure}

\paragraph{Convex envelope.}
The \emph{convex envelope} of $f$ is
\[
f^{**}(x) = \sup \{ g(x) \mid g \le f, \; \text{$g$ convex} \},
\]
\ie, its value at a point $x$ is the greatest value of $g(x)$
obtained for any convex function $g$ that minorizes $f$.
The convex envelope is also the (Fenchel) conjugate of the conjugate of $f$, \ie, $(f^*)^*$, 
where the superscript $*$ denotes conjugation \cite[\S 5]{rockafellar1970convex}.
(This explains why we denote the convex envelope of $f$ as $f^{**}$.)
If $f$ is convex, then we have $f^{**} = f$.
An example is shown in figure~\ref{f-conv-env-example}.

\paragraph{Bound from convex relaxation.}
By replacing the separable functions $f_i$ in problem~\eqref{e-sap}
with their convex envelopes $f_i^{**}$,
we obtain the problem
\begin{equation}
  \label{e-sap-relax}
  \begin{array}{ll}
    \mbox{minimize} & \sum_{i=1}^n f_i^{**}(x_i)\\
    \mbox{subject to} & Ax = b,
  \end{array}
\end{equation}
which we call the \emph{convex relaxation} of \eqref{e-sap}.
The relaxed problem~\eqref{e-sap-relax} is a convex separable--affine problem.
The optimal value of this relaxation is denoted $d^\star$.
(The reason for this notation will become clear later.)
Because $f_i^{**} \le f_i$, we have $d^\star \le p^\star$.
Indeed, if $f_i$ is convex for all $i$,
then $f_i = f_i^{**}$ and therefore $d^\star = p^\star$.

\paragraph{Dual problem.}
The dual problem of \eqref{e-sap} is
\begin{equation}
  \label{e-dual-prob}
  \begin{array}{ll}
    \mbox{maximize} & \lambda^T b -\sum_{i=1}^n f_i^*(-\nu_i) \\
    \mbox{subject to}
        & A^T \lambda = \nu.
  \end{array}
\end{equation}
The variables are $\nu\in\reals^n$ and $\lambda\in\reals^m$.
Because the linear term in \eqref{e-dual-prob} is separable,
the dual problem is itself a separable--affine problem.

The dual problem is always a convex optimization problem,
even when the primal problem~\eqref{e-sap} is not.
The optimal value of the dual problem is $d^\star$,
which is the optimal value of the relaxed problem~\eqref{e-sap-relax}.
Weak duality, which always holds, states that $d^\star \le p^\star$.
If \eqref{e-sap} is convex,
then strong duality holds, which means that $d^\star = p^\star$
\citep[Prop.\ 5.2.1]{bertsekasnonlinear}.

\section{ADMM}
\label{s-admm}
We can apply the alternating direction method of multipliers (ADMM)
to the separable--affine problem~\eqref{e-sap}.

\paragraph{ADMM-form problem.}
We start with the equivalent problem
\begin{align}
\label{e-admm-prob}
\begin{array}{ll}
\text{minimize}
    & \mathcal{I}_\mathcal{A}(z) + \sum_{i=1}^n f_i(x_i) \\
\text{subject to}
    & x = z
\end{array}
\end{align}
with variables $x\in\reals^n$ and $z\in\reals^n$.
Here $\mathcal I_{\mathcal A}$ is the indicator function over the affine constraints of (\ref{e-sap}), \ie,
\[
\mathcal{I_{\mathcal A}}(z) = \begin{cases}
0 & \text{if $Az = b$} \\
\infty & \text{otherwise.}
\end{cases}
\]

\paragraph{Augmented Lagrangian.}
The augmented Lagrangian of (\ref{e-admm-prob}) is
\[
L(x, z, \lambda) = f(x) + I_{\mathcal A}(z) + \frac12 \| x - z + \lambda \|^2.
\]
Our definition of the augmented Lagrangian
lacks the parameter $\rho>0$ from the standard definition given by \citet{bertsekasnonlinear}.
Including this parameter is equivalent to using problem scaling parameter $D = \rho I$,
as discussed in section~\ref{s-prob-scaling}.

\paragraph{ADMM iterations.}
\label{methods:admm:iterations}
The ADMM algorithm iterates $x^k$, $z^k$, and $\lambda^k$, for $k = 0, 1, 2, \dots$, are
\begin{gather}
\label{e-admm-iterates}
\begin{aligned}
  x^{k+1} &= \argmin_{x} L(x, z^k, \lambda^k) \\
  z^{k+1} &= \argmin_{z} L(x^{k+1}, z, \lambda^k) \\
  \lambda^{k+1} &= \lambda^k + x^{k+1} - z^{k+1}.
\end{aligned}
\end{gather}
The initial values are $z^0$ and $\lambda^0$.

\paragraph{Convergence.}
\label{s-convergence}
If $f$ is convex and a solution to \eqref{e-sap} exists,
then as $k\to\infty$, we have $f(x^k) \to p^\star$, $Ax^k \to b$,
and $\lambda^k \to \lambda^\star$,
where $\lambda^\star$ is an optimal dual variable to \eqref{e-sap}.
In the general case when $f$ is nonconvex, there is no such guarantee.

\paragraph{Solving the linear system.}
Here we describe the $z^{k+1}$ update in (\ref{e-admm-iterates}).
Minimizing the augmented Lagrangian involves solving the equality-constrained
least-squares problem
\begin{align*}
\begin{array}{ll}
\text{minimize} & \| z - x^{k+1} + \lambda^k \|^2 \\
\text{subject to} & Az = b
\end{array}
\end{align*}
with variable $z$.
The minimizing $z$ (which becomes $z^{k+1}$)
can be found by solving the linear system of equations
\begin{equation}
\begin{bmatrix} I & A^T \\ A & 0 \end{bmatrix}
\begin{bmatrix} z^{k+1} \\ \nu \end{bmatrix}
=
\begin{bmatrix} x^k - \lambda^k \\ b \end{bmatrix}.
\label{e-linear-system}
\end{equation}
For each iteration in the ADMM algorithm,
we solve \eqref{e-linear-system}
for different values of the the right-hand side
(\ie, for different values of $x^k$ and $\lambda^k$.)
We can do this efficiently by factorizing this matrix before the first iteration
and caching it for repeated use \citep[\S 3.1]{osqp}.

\paragraph{Separable update.}
\label{s-separable-update}
The update rule for $x^{k+1}$ in \eqref{e-admm-iterates} can be written
\[
x^{k+1} = \argmin_{x} \left( f(x) + \frac12 \| x - z^k + \lambda^k \|^2
\right).
\]
Because $f$ is separable,
we can perform this update by solving $n$ univariate optimization problems. 
The update for each component $x_i$ is given by
\[
x_i^{k+1} = \argmin_{x_i} \left( f_i(x_i) + \frac12
\left(x_i - z_i^k + \lambda_i^k \right)^2 \right).
\]
These updates can be expressed using the proximal operator of $f_i$:
\[
x_i^{k+1} = \prox_{f_i}\big(z_i^k - \lambda_i^k\big).
\]
(For more information, see \cite{boyd2011distributed}.)
Several methods can be used to solve these small problems,
including exhaustive search.
For the portfolio problem~\eqref{e-port-opt},
$f_i$ is piecewise quadratic,
and the proximal operator can be readily computed, as shown
in section~\ref{e-pwq-prox}.

\paragraph{Initialization.}
\label{s-admm-initialization}
If $f$ is convex, the ADMM iterates converge to a solution of \eqref{e-sap}
regardless of the initialization.
However, if $f$ is not convex,
the choice of initialization can have a large impact on the result of ADMM.

In this case, one good choice is to initialize the iterates
using the solution to the relaxed problem \eqref{e-sap-relax}.
This has the added benefit of providing a lower bound on the optimal problem value,
which can be used to judge the quality of the iterates produced by ADMM.

\paragraph{Scaling.}
Although the scaling parameters $D$ and $E$
do not change the solution of \eqref{e-sap},
they can have a large impact on the rate of convergence of ADMM.
(When \eqref{e-sap} is not convex, they can also affect the quality of the iterates.)
One effective method to choose $D$ and $E$ is through a simple hyperparameter search,
in which we and tune $D$ and $E$ to minimize the run time over a sampling of similar problems.
This requires expressing $D$ and $E$ in terms of a small number of free 
parameters, as in done for example in section~\ref{s-prob-data}.
Another approach is to choose $D$ and $E$ to reduce the condition number of $A$,
\eg, through equilibration \citep[\S 5]{osqp}.

\paragraph{Termination criteria.}
We need termination criteria that work well in the case when $f$ is not convex.
Termination criteria for ADMM, when applied to convex problems,
are discussed in \cite[\S 3.3]{boyd2011distributed}.

For any candidate point $x \in \reals^n$ satisfying $Ax = b$,
we define the pseudo-objective and residual values as
\begin{align*}
    o(x) &= f\big(\Pi_{\dom(f)}(x)\big) \\
    r(x) &= \dist(x, \dom(f)).
\end{align*}
We use the following procedure to check for termination:
\begin{enumerate}
\item From the current iterates, find a point $x$ satisfying $Ax = b$. \label{s-step-1}
\item Compute $o(x)$ and $r(x)$.
\item If $r(x) < \epsilon^{\rm res}$ and $o(x) < o^{\rm best}$, then update $o^{\rm best}$
    and the best iterates.
\item Terminate if the best objective $o^{\rm best}$ has not improved by more than $\epsilon^{\rm obj}$
    in more than $N$ iterations.
\end{enumerate}
There are several ways to carry out step 1.
For example, we can take the iterate $x^k$,
which satisfies $Ax^k = b$.
For the portfolio optimization problem~\eqref{e-port-opt},
we can take the elements of the iterate $z^k$ that represent the asset holdings $h$
and compute the corresponding cash amount $c$ and factor exposures $y$.
Then the vector $x = (h, c, y)$ satisfies $Ax = b$.

\section{Implementation}
\label{s-implementation}
Carrying out ADMM requires evaluating the proximal operator of univariate functions,
and computing the bound from \ref{s-bounds} requires computing the convex envelope.
In this section, we show how to carry out these two operations
for the specific case in which the univariate functions are piecewise quadratic.
We use the notation $\varphi$ for a generic univariate function,
which stands in for $f_i$ in \eqref{e-sap}.

\subsection{Piecewise quadratic functions}
\label{impl:pwq}
Consider the piecewise quadratic function $\varphi$ with $k$ pieces, defined as
\begin{equation}\label{e-pwq}
  \varphi(x) = \begin{cases}
    p_1x^2 + q_1x + r_1 & \text{if $x \in [a_1, b_1]$} \\
    \hfill\vdots\hfill & \hfill\vdots\hfill \\
    p_kx^2 + q_kx + r_k & \text{if $x \in [a_k, b_k]$} \\
    +\infty & \text{otherwise} \\
  \end{cases}
\end{equation}
with $a_1 \le b_1 \le a_2 \le \cdots \le a_k \le b_k$.
We denote by $\varphi_i$ the $i$th piece of $\varphi$, \ie,
\begin{equation}\label{e-pwq-piece-def}
\varphi_i(x) =
\begin{cases}
p_ix^2 + q_ix + r_i & \text{if $x \in [a_i, b_i]$} \\
+\infty & \text{otherwise.} \\
\end{cases}
\end{equation}
Because we do not require $\varphi$ to be continuous,
and each piece is defined over a closed interval,
$\varphi$ may be multiply defined for some values of $x$ (at the boundaries of some
intervals).
To remedy this, we use the convention that $\varphi(x)$
is the minimum over all such candidate values, \ie,
$\varphi(x) = \min \{ \varphi_i(x) \mid \text{$\forall i$ s.t.\ $x \in [a_i, b_i]$} \}$.
We also assume the description of $\varphi$ is irreducible,
\ie, none of the $k$ pieces can be dropped without changing the function value
at some point.

\subsection{Proximal operator}
\label{e-pwq-prox}
The proximal operator of $\varphi$ is
\begin{equation}
\prox_{\varphi}(u) = \argmin_x \left(\varphi(x) + \frac12 (x - u)^2 \right).
\end{equation}
Evaluating the proximal operator is done in two steps.
First, we compute the piecewise quadratic function
$\varphi(x) + (x - u)^2/2$,
which is done by adding $1/2$, $-u$, and $u^2/2$
to the coefficients $p_i$, $q_i$, and $r_i$, for all $i$.
Computing the minimizer of this function can
be done by computing the minimum value of each piece,
taking the minimum over these values,
and then finding a value of $x$ that attains this minimum.

If $\varphi$ is convex, its proximal operator of $\varphi$ can be expressed as
\begin{align*}
& \prox_\varphi(u) = \hfill \\
    & \qquad \begin{cases}
      a_1 & \text{if $u \in (-\infty, (2p_1 + 1)a_1 + q_1]$}\\
      (u - q_1)/(1 + 2p_1) & \text{if $u \in [(2p_1 + 1)a_1 + q_1, (2p_1 + 1)b_1 + q_1]$}\\
      \hfill\vdots\hfill & \hfill\vdots\hfill\\
      a_j & \text{if $u \in [(2p_{j-1} + 1)b_{j-1} + q_{j-1}, (2p_j + 1)a_j + q_j]$}\\
      (u - q_j)/(1 + 2p_j) & \text{if $u \in [(2p_j + 1)a_j + q_j, (2p_j + 1)b_j + q_j]$}\\
      \hfill\vdots\hfill & \hfill\vdots\hfill\\
      a_k & \text{if $u \in [(2p_{k-1} + 1)b_{k-1} + q_{k-1}, (2p_k + 1)a_k + q_k]$}\\
      (u - q_k)/(1 + 2p_k) & \text{if $u \in [(2p_k + 1)a_k + q_k, (2p_k + 1)b_k + q_k]$}\\
      b_k & \text{if $u \in [(2p_k + 1)b_k + q_k, \infty)$}.
    \end{cases}
\end{align*}
We note that some of these intervals may be degenerate.

\subsection{Convex envelope}
\label{s-cvx-env-pwq}
\begin{figure}
\centering
\ifmakeplots
    \begin{tikzpicture}

        \begin{groupplot}[group style={group size=4 by 1, horizontal sep=0.0266\textwidth},
                          xmin=-1,
                          xmax=8,
                          ymin=-4,
                          width = 0.23\textwidth,
                          height = 0.23\textwidth / \aspectratio,
                          scale only axis=true,
                          ymax=5,
                          %y ticks=none,
                          ytick=\empty,
                          legend cell align={left},
                          ]

        \nextgroupplot[
            legend to name={CommonLegend2},
            legend style={legend columns=4, /tikz/every even column/.append style={column sep=0.5cm}},
            legend cell align={left},
            xtick = {3,4},
            title=\text{$i=1$},
        ]
            \addplot [thick, black, dotted] table{\datadir/example_function.csv};
            \addplot [thick, blue] table[x=X, y=Y, col sep=comma]{\datadir/recursive_envelope_glast_2.csv};
            %\addplot [thick, dotted, red] table[x=X, y=Y, col sep=comma]{\datadir/recursive_envelope_fi_2.csv};
            \addplot [thick, dashed, blue] table[x=X, y=Y, col sep=comma]{\datadir/recursive_envelope_gnext_2.csv};

            \addlegendentry{$\varphi(x)$}
            \addlegendentry{$\psi^i(x)$}
            \addlegendentry{$\psi^{i+1}(x)$}

        \nextgroupplot[
            title=\text{$i=2$},
            xtick = {4,6},
        ]
            \addplot [thick, black, dotted] table{\datadir/example_function.csv};
            \addplot [thick, blue] table[x=X, y=Y, col sep=comma]{\datadir/recursive_envelope_glast_3.csv};
            %\addplot [thick, dotted, red] table[x=X, y=Y, col sep=comma]{\datadir/recursive_envelope_fi_3.csv};
            \addplot [thick, dashed, blue] table[x=X, y=Y, col sep=comma]{\datadir/recursive_envelope_gnext_3.csv};

        \nextgroupplot[
            title=\text{$i=3$},
            xtick = {6, 7.5},
        ]
            \addplot [thick, black, dotted] table{\datadir/example_function.csv};
            \addplot [thick, blue] table[x=X, y=Y, col sep=comma]{\datadir/recursive_envelope_glast_4.csv};
            %\addplot [thick, dotted, red] table[x=X, y=Y, col sep=comma]{\datadir/recursive_envelope_fi_4.csv};
            \addplot [thick, dashed, blue] table[x=X, y=Y, col sep=comma]{\datadir/recursive_envelope_gnext_4.csv};

        \nextgroupplot[
            title=\text{$i=4$},
            xtick = {7.5},
        ]
            \addplot [thick, black, dotted] table{\datadir/example_function.csv};
            \addplot [thick, blue] table[x=X, y=Y, col sep=comma]{\datadir/recursive_envelope_glast_5.csv};
            %\addplot [thick, dotted, red] table[x=X, y=Y, col sep=comma]{\datadir/recursive_envelope_fi_5.csv};
            \addplot [thick, dashed, blue] table[x=X, y=Y, col sep=comma]{\datadir/recursive_envelope_gnext_5.csv};

            %\legend{$\psi^2(x)$, $\varphi_3(x)$, $\psi^3(x)$}

        \end{groupplot}
        \path ($(group c1r1.south) - (0cm, 0.8cm)$) -- node[below]{\ref{CommonLegend2}} ($(group c4r1.south) - (0cm, 0.8cm)$);

    \end{tikzpicture}
\fi
\caption{
Recursive computation of the convex envelope for the function shown in figure~\ref{f-conv-env-example}.
}
\label{f-recursive-envelope}
\end{figure}
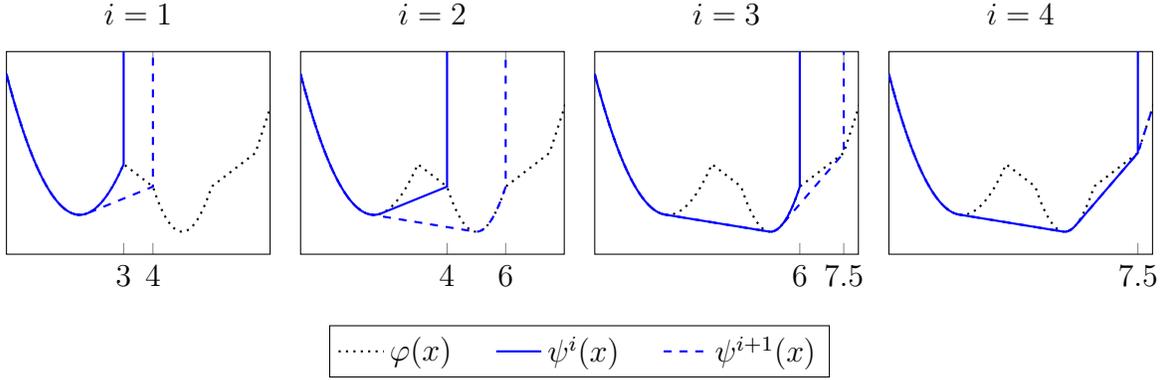

In this section we show how to compute the convex envelope of $\varphi$,
which is required for computing the bound discussed in section~\ref{s-bounds}.
This follows the same lines as \cite{gardiner2010convex}.
To do this, first note that the convex envelope $\varphi^{**}$ can be computed recursively:
\begin{gather}
\label{e-envelope-recursion}
\begin{aligned}
\varphi^{**} &= \min\{\psi^i, \varphi_i\}^{**} \\
\psi^i &= \min\{\psi^{i-1}, \varphi_{i-1}\}^{**} \\
& \qquad \vdots \\
\psi^2 &= \min\{\psi^1, \varphi_1\}^{**} \\
\psi^1 &= \varphi_1.
\end{aligned}
\end{gather}
Here $\min$ is the pointwise minimum operation between functions.
An example of this recursion for the function
\begin{equation}
\label{e-pwq-example}
  \varphi(x) = \begin{cases}
    x^2 - 3x - 3 & \text{if $x \in [-\infty, 3]$}\\
    -x + 3 & \text{if $x \in [3, 4]$}\\
    2x^2 - 20x + 47 & \text{if $x \in [4, 6]$}\\
    x - 7 x & \text{if $\in [6, 7.5]$}\\
    4x + 29 & \text{if $x \in [7.5, \infty]$},
  \end{cases}
\end{equation}
is shown in figure~\ref{f-recursive-envelope}.

\begin{figure}
\centering
\ifmakeplots
    \begin{tikzpicture}
        \begin{axis}[xmin=-1.5,
                     xmax=7.5,
                     ymin=-6,
                     ymax=3,
                     width = 0.45\textwidth,
                     height = 0.45\textwidth / \aspectratio,
                     scale only axis=true,
                     ytick=\empty,
                     xtick = {1.39366479822405220368, 4.94683239911202576877, 6},
                     xticklabels = {$x_\varphi$, $x_\psi$, $b_2$},
                     legend cell align={left}
                     ]

            \addplot [thick, black, dotted] table{\datadir/example_function.csv};
            \addplot [thick, blue] table[x=X, y=Y, col sep=comma]{\datadir/recursive_envelope_glast_3.csv};
            %\addplot [thick, dotted, red] table[x=X, y=Y, col sep=comma]{\datadir/recursive_envelope_fi_3.csv};
            \addplot [thick, dashed, blue] table[x=X, y=Y, col sep=comma]{\datadir/recursive_envelope_gnext_3.csv};

            \coordinate (brace1l) at (axis cs:-1.5, -2.5); % {brace1l};
            \coordinate (brace1r) at (axis cs:1.393664, -2.5); % {brace1r};
            \coordinate (brace2r) at (axis cs:4.946832, -2.5); % {brace1r};
            \coordinate (brace3r) at (axis cs:6, -2.5); % {brace1r};

            \draw [decorate,decoration={brace,amplitude=10pt,mirror,raise=6mm}]
            (brace1l) -- (brace1r) node [black,midway, below, yshift=-10mm] {$\psi^i$};
            \draw [decorate,decoration={brace,amplitude=10pt,mirror,raise=6mm}]
            (brace1r) -- (brace2r) node [black,midway, below, yshift=-10mm] {$\alpha x + \beta$};
            \draw [decorate,decoration={brace,amplitude=10pt,mirror,raise=6mm}]
            (brace2r) -- (brace3r) node [black,midway, below, yshift=-10mm] {$\varphi_i$};

            %\addplot [thick, red] table[x=X, y=Y, col sep=comma]{\datadir/recursive_envelope_gnext_left.csv};
            %\addplot [thick, yellow] table[x=X, y=Y, col sep=comma]{\datadir/recursive_envelope_gnext_middle.csv};
            %\addplot [thick, magenta] table[x=X, y=Y, col sep=comma]{\datadir/recursive_envelope_gnext_right.csv};

        \end{axis}
    \end{tikzpicture}
\fi
\caption{
The three components that make up $\psi^{i+1}$ according to equation~\eqref{e-env-param-cond-pwq},
with $\varphi$ given in \eqref{e-pwq-example} and $i=2$.
}
\label{f-conv-env-components}
\end{figure}
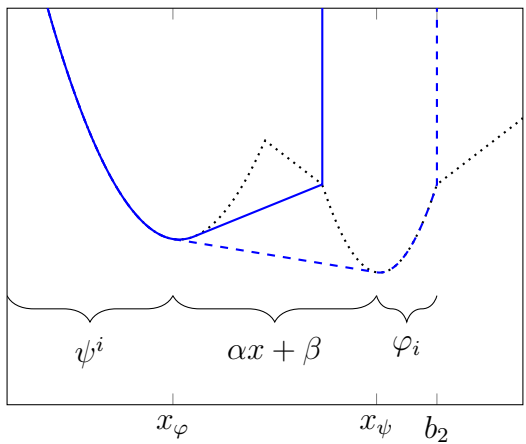

We now discuss how to
carry out each line in the recursion~\eqref{e-envelope-recursion},
\ie, how to compute $\psi^{i+1}$ given $\psi^i$.
Note that $\psi^i$ and $\varphi_i$ are convex functions.
The graph of $\psi^i$ lies to the left of the graph of $\varphi_i$,
meaning that $x \le z$ for all $x\in \dom \psi^i$ and $z \in \dom \varphi_i$.
Therefore $\psi^{i+1}$ has the simple form
\begin{equation}
  \label{e-env-two-cvx-funs}
  \psi^{i+1}(x)
  = \begin{cases}
    \psi^i(x) & \text{if $x \in [a_1, x_\psi]$} \\
    \alpha x + \beta & \text{if $x \in [x_\psi, x_\varphi]$} \\
    \varphi_i(x) & \text{if $x \in [x_\varphi, b_i]$}
  \end{cases}
\end{equation}
for some $\alpha \in \reals$, $\beta \in \reals$,
$x_\varphi \in \dom \varphi_i$, and $x_\psi \in \dom \psi^i$.
This is shown in figure~\ref{f-conv-env-components}.
We allow for the case when $x_\varphi = -\infty$ or $x_\psi = \infty$;
in these cases, the first or last interval is degenerate, and can be ignored.
The parameters $\alpha$ and $\beta$ are unique;
$x_\varphi$ and $x_\psi$ need not be.
Finding these parameters is straightforward but tedious;
further details are given in appendix~\ref{s-envelope-appendix}.

\section{Numerical example}
\label{s-example}
In this section we provide an example of the method applied to tax-aware portfolio optimization with small account sizes.
We use a passive, index-tracking strategy with $\alpha = 2\gamma^{\rm risk} Vh^{\rm bm}$,
where the elements of $h^{\rm bm} \in\reals^l$
are the weights of a benchmark portfolio.
This means that minimizing the objective term $\alpha^T h + \gamma^{\rm risk} h^T V h$
is equivalent to minimizing the (squared) active risk
$(h - h^{\rm bm})^T V (h - h^{\rm bm})$.

The trading cost combines the bid--ask spread model
with a penalty on the number of trades and capital gains:
\[
\phi^{\rm trd}(u) = \gamma^{\rm sprd} s^T |u| + c^{\rm trd}\card(u) + \gamma^{\rm tax} L(u).
\]

The holding cost combines the asset position limits~\eqref{e-position-limits}
and the per-asset holding cost~\eqref{e-holding-cost}.
We use the lower bound $h^{\rm lb} = 0$, which encodes a long-only constraint,
and the upper bound for asset $i$ is
$h^{\rm ub}_i = \max \{ 3 h^{\rm bm}_i, h^{\rm init}_i \}$,
\ie, it is the greater of the current holdings and 3 times the benchmark weight.
The combined holding cost is
\[
\phi^{\rm hld}(h) = \begin{cases}
c^{\rm hld}\card(h)  & \text{if $0 \le h \le h^{\rm ub}$} \\
\infty & \text{otherwise.} \\
\end{cases}
\]

\subsection{Problem data}
\label{s-prob-data}

\paragraph{Backtest setup.}
We generate instances of problem~\eqref{e-port-opt}
from backtests of a tax loss harvesting strategy.
Our backtest dataset consists of 204 months over a 17 year period from August 2002 through August 2019.
We use this dataset to carry out 12 staggered six-year-long backtests,
with the first starting in August 2002 and ending in July 2008,
and the the last starting in August 2013 and ends in July 2019.
Every month, a single instance of problem~\eqref{e-port-opt} is solved to rebalance the portfolio.
We then save this problem instance, which we use to evaluate our heuristic and bound.
The full details of this setup are described by \citet[\S 6]{moehle2020tax},
including the sources of data used and the exact timing of the rebalance trades.
(The main difference between our formulation and theirs
is the addition of the per-asset holding and trading cost terms.)

When solving problem~\eqref{e-port-opt}, we use the parameters
\[
\gamma^{\rm risk} = 100, \qquad c^{\rm hld} =  c^{\rm trd} = 3\times 10^{-5},
\qquad
\gamma^{\rm tax} = 1.
\]
The parameters $\gamma^{\rm trd}$ and $\gamma^{\rm hld}$ are redundant, 
and were set to 1.
With these values, the portfolio maintains active risk around 0.5\%--1\%,
holds only 200--300 securities of the S\&P 500,
and typically trades around 30 securities per month.
The invested fraction is maintained between
$\eta^{\rm lb} = 0.98$ and $\eta^{\rm ub} = 0.99$.

\paragraph{Generated problems.}
The procedure given above resulted in 692 instances of problem~\eqref{e-port-opt}.
The mean optimal utility $U^\star$ across these problems
ranged from $-1028$ to $184$ basis points ($0.0001$, one hundredth of one percent)
with a mean and standard deviation of 
$94$ basis points and $57$ basis points, respectively.

\paragraph{Algorithm parameters.}
We used the scaling parameters 
\[
D = \diag(100 \ones_l, 3, 100 \ones_k), \qquad
E = \diag(100 \ones_m),
\]
where $\ones_l$ and $\ones_k$ are the $l$ and $k$ dimensional vectors with 
all entries one.
We used the stopping criterion parameters 
\[
\epsilon^{\rm res} = 3 \times 10^{-4}, \qquad \epsilon^{\rm obj} = 10^{-5}.
\]
We check the termination conditions once every 10 steps,
and terminate if the objective have not improved in more than $N = 50$ iterations.

\subsection{Results}
\label{s-results}

\paragraph{Heuristic quality.}
All of the 692 problems we ran the algorithm on converged.
To evaluate the ADMM heuristic,
we compare the objective values obtained by ADMM on the 692 problem instances,
which we denote $p_{\rm admm}$,
to the lower bounds $d^\star$
obtained by solving the convex relaxation~\eqref{e-sap-relax}.
The optimal problem value $p^\star$ lies between these values, \ie,
$$
d^\star \le p^\star \le p_{\rm admm}.
$$

\begin{figure}
\centering
\ifmakeplots
    \begin{tikzpicture}
        \begin{axis}[xlabel=$p_{\rm admm} - d^\star$,
                     ylabel=instances,
                     width = 0.5\textwidth,
                     height = 0.4\textwidth / \aspectratio,
                     scale only axis=true,
                     legend cell align={left},
                     area style,
                     ymin = 0,
                     xticklabels = {$10^{-7}$, $10^{-6}$, $10^{-5}$, $10^{-4}$, $10^{-3}$},
                     xtick = { -7, -6, -5, -4, -3 },
                     ]

            \addplot+[ybar interval, mark=no, draw opacity=0.5, fill opacity=0.4, fill = blue]
                      table[x=diff, y=freq, col sep=comma]{\datadir/gaps.csv};
        \end{axis}
    \end{tikzpicture}
\fi
\caption{Distribution of suboptimality gaps $p_{\rm admm} - d^\star$.}
\label{f-gaps}
\end{figure}
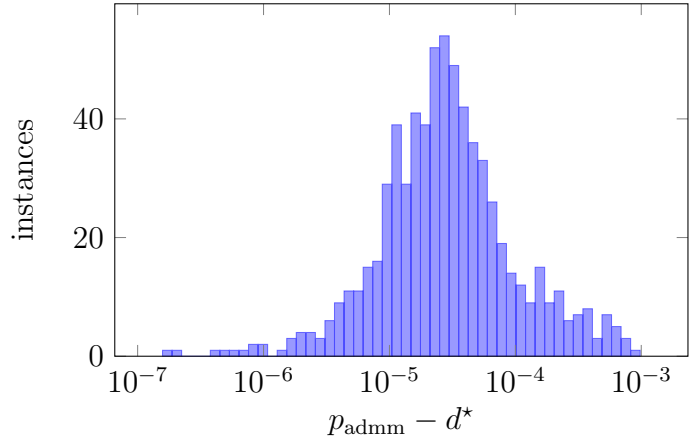

Figure~\ref{f-gaps} shows the differences $p_{\rm admm} - d^\star$.
These range from $0$ to $10$ basis points, with mean $0.6$ basis points and
standard deviation $1.1$ basis points.
These values are quite small compared to the problems values $p^\star$, which range from
$-1028$ to $184$ basis points.
This implies that the ADMM heuristic produces nearly optimal points on all 692 problem instances.

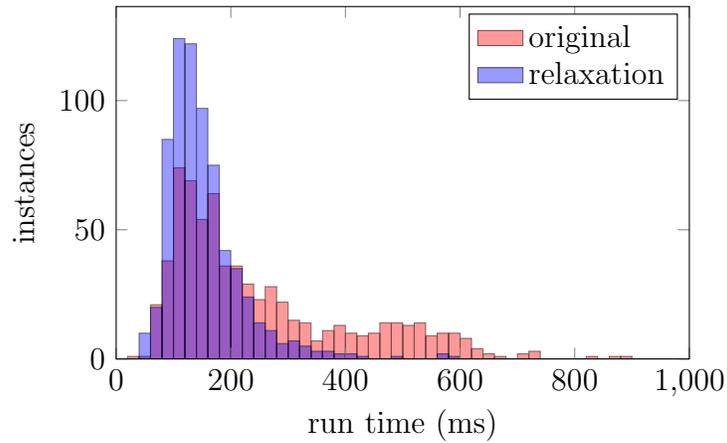
\begin{figure}
\centering
\ifmakeplots
    \begin{tikzpicture}
        \begin{axis}[xlabel=run time (ms),
                     ylabel=instances,
                     width = 0.5\textwidth,
                     height = 0.4\textwidth / \aspectratio,
                     scale only axis=true,
                     legend cell align={left},
                     area style,
                     ymin=0,
                     xmin=0,
                     xmax=1000,
                     ]

            \addplot [
                hist={data=x, bins=50, data min=0, data max=1000},
                fill = red,
                fill opacity=0.4,
                draw opacity=0.5,
            ]
            table[x=uppertime, col sep=comma] {\datadir/runtimes_raw.csv};

            \addplot [
                hist={data=x, bins=50, data min=0, data max=1000},
                fill = blue,
                fill opacity=0.4,
                draw opacity=0.5,
            ]
            table[x=lowertime, col sep=comma] {\datadir/runtimes_raw.csv};

            \legend{original, relaxation}
        \end{axis}
    \end{tikzpicture}
\fi
    \caption{
        Run time distributions for the ADMM algorthm applied to the
        original problem~\eqref{e-port-opt} as well as the relaxed problem~\eqref{e-sap-relax}.
    }
\label{f-runtimes}
\end{figure}

\paragraph{Speed.} In addition to being accurate, ADMM was also fast. 
It solved the original (nonconvex)
problems in an average of $251$ ms (with standard deviation $159$ ms),
and solved the relaxed (convex) problems
in an average of $152$ ms (with a standard deviation $67$ ms).
Figure~\ref{f-runtimes} shows the distributions of run times
for original and relaxed problems.

\subsection{Single problem instance}
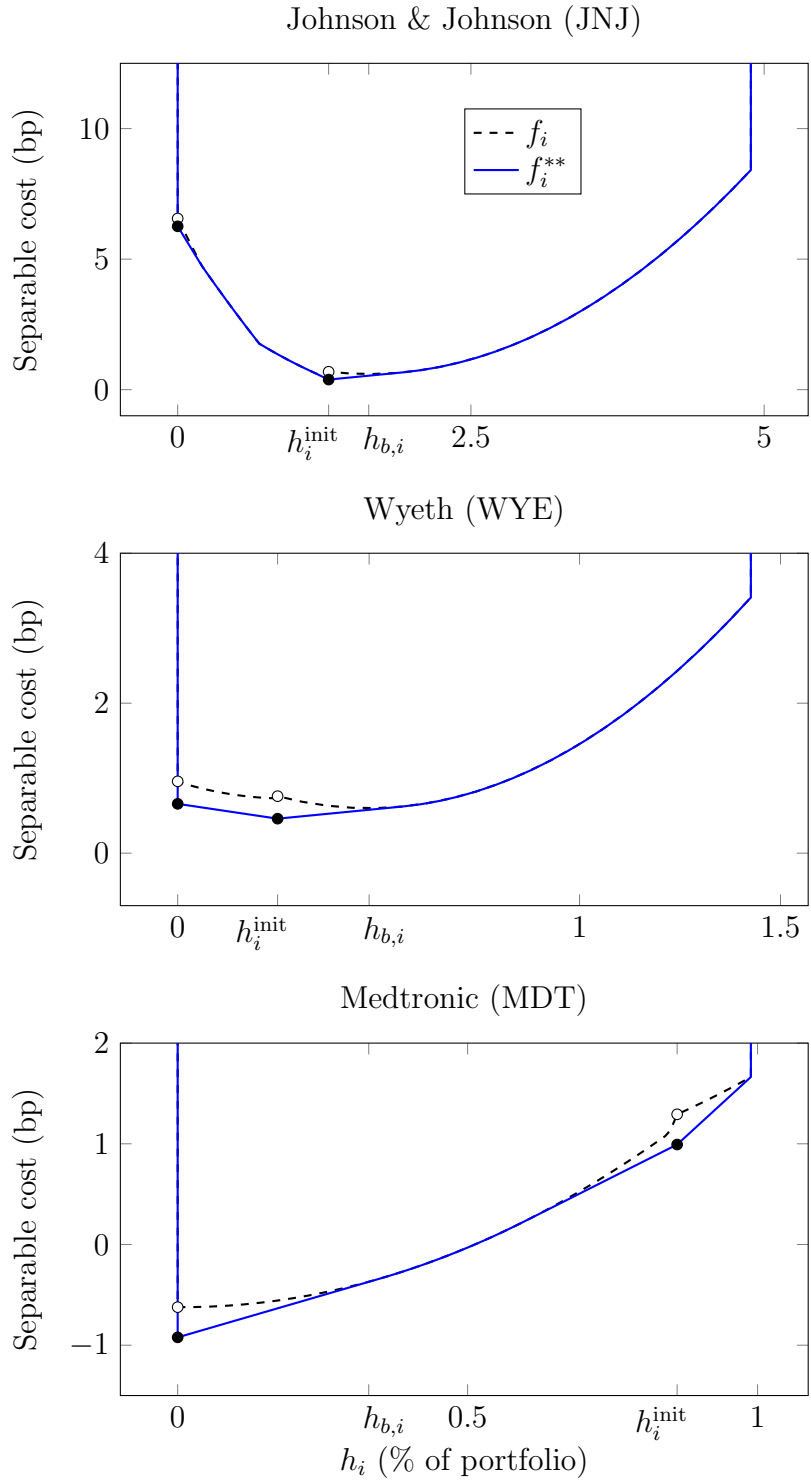
\begin{figure}
\centering
\ifmakeplots
    \begin{tikzpicture}
        \begin{groupplot}[group style={group size=1 by 3, vertical sep=0.12\textwidth},
                          width = 0.6\textwidth,
                          height = 0.4\textwidth / \aspectratio,
                          scale only axis=true,
                          legend cell align={left},
                          legend style={at={(0.5,0.75)},anchor=west},
                          ylabel={Separable cost (bp)},
                          %xmin=-.5,
                          %xmax=5.5,
                          %ytick scale label code/.code={$x$ $ ( 10^{#1}$ m)},
                          y filter/.code= {\pgfmathmultiply{#1}{10000}},
                          x filter/.code= {\pgfmathmultiply{#1}{100}}
                         ]
                \nextgroupplot[
                    ymin=-1,
                    ymax=12.5,
                    title={Johnson \& Johnson (JNJ)},
                    xtick = {0, 1.28692, 1.62891, 2.5, 5},
                    xticklabels = {0, $h^{\rm init}_i \quad$, $\quad h_{b, i}$, 2.5, 5}
                ]
                    \addplot [only marks, color=black, forget plot, mark options={fill=white}, on layer=background1] coordinates {(0.012869220317084354,6.860903948553906e-05)};
                    \addplot [only marks, color=black, forget plot, mark options={fill=white}, on layer=background1] coordinates {( 0, 0.0006556585832405394)};
                    \addplot [only marks, color=black, forget plot, on layer=background2] table[x=px, y=py, col sep=comma]{\datadir/orig_points_input23_asset467.csv};
                    \addplot [thick, black, dashed] table[x=lx, y=ly, col sep=comma]{\datadir/orig_lines_input23_asset467.csv};
                    \addplot [thick, blue] table[x=lx, y=ly, col sep=comma]{\datadir/env_input23_asset467.csv};
                    %\addplot +[mark=none, black] coordinates {(0.0128692, -0.0001) (0.0128692, 0.002)};
                    %\addplot +[mark=none, black, dashed] coordinates {(0.0162891, -0.0001) (0.0162891, 0.002)};
                    \legend{$f_i$, $f_i^{**}$, $h^{\rm init}_i$, $h_{b, i}$};

                \nextgroupplot[
                    ymin=-.7,
                    ymax=4,
                    xtick = {0, .24861, .4764, 1.0, 1.5},
                    xticklabels = {0, $h^{\rm init}_i\quad$, $\quad h_{b, i}$, 1, 1.5},
                    title={Wyeth (WYE)},
                    %\addplot +[mark=none, black] coordinates {(0.0024861, -0.0001) (0.0024861, 0.0005)};
                    %\addplot +[mark=none, black, dashed] coordinates {(0.004754, -0.0001) (0.004754, 0.0005)};
                 ]
                    \addplot [only marks, color=black, forget plot] table[x=px, y=py, col sep=comma, on layer=background2]{\datadir/orig_points_input3_asset120.csv};
                    \addplot [only marks, color=black, forget plot, mark options={fill=white}, on layer=background1] coordinates {(0.0024861058023495537, 7.598589033128831e-05)};
                    \addplot [only marks, color=black, forget plot, mark options={fill=white}, on layer=background1] coordinates {( 0, 9.578308360860615e-05)};
                    \addplot [thick, black, dashed] table[x=lx, y=ly, col sep=comma]{\datadir/orig_lines_input3_asset120.csv};
                    \addplot [thick, blue] table[x=lx, y=ly, col sep=comma]{\datadir/env_input3_asset120.csv};
                    %\addplot +[mark=none, black] coordinates {(0.0024861, -0.0001) (0.0024861, 0.0005)};
                    %\addplot +[mark=none, black, dashed] coordinates {(0.004754, -0.0001) (0.004754, 0.0005)};

                \nextgroupplot[
                    ymin=-1.5,
                    ymax=2,
                    title={Medtronic (MDT)},
                    xlabel={$h_i$ (\% of portfolio)},
                    xtick = {0, 0.32949, 0.5, 0.86131, 1.0, 1.5},
                    xticklabels = {0, $\quad h_{b, i}$, 0.5, $h^{\rm init}_i\quad$, 1, 1.5}
                ]
                    \addplot [only marks, color=black, forget plot] table[x=px, y=py, col sep=comma, on layer=background2]{\datadir/orig_points_input22_asset564.csv};
                    \addplot [only marks, color=black, forget plot, mark options={fill=white}, on layer=background1] coordinates {(0.008613192324435248,0.00012927195554923655)};
                    \addplot [only marks, color=black, forget plot, mark options={fill=white}, on layer=background1] coordinates {( 0, -6.221843946052044e-05)};
                    \addplot [thick, black, dashed] table[x=lx, y=ly, col sep=comma]{\datadir/orig_lines_input22_asset564.csv};
                    \addplot [thick, blue] table[x=lx, y=ly, col sep=comma]{\datadir/env_input22_asset564.csv};
                    %\addplot +[mark=none, black] coordinates {(0.008613, -0.00012) (0.008613, 0.0003)};
                    %\addplot +[mark=none, black, dashed] coordinates {(0.0032949, -0.00012) (0.0032949, 0.0003)};
        \end{groupplot}
    \end{tikzpicture}
\fi
\caption{
The cost functions $f_i$ and their relaxations $f_i^{**}$
for the stocks of Johnson \& Johnson, Wyeth, and Medtronic,
for a single simulated trade.}
\label{f-sec-fi}
\end{figure}

Figure~\ref{f-sec-fi} shows the separable cost functions $f_i$ corresponding to
the stocks of Johnson \& Johnson (JNJ), Wyeth (WYE), and Medtronic (MDT).
We show the initial holdings $h^{\rm init}_i$,
the benchmark holdings $h^{\rm bm}_i$, and the optimal post-trade holdings $h_i^\star$.

The upward curvature of all three functions is primarily due to the specific risk term
$D_{ii} (h_i - h^{\rm bm}_i)^2$.
Each $f_i$ has two discontinuities:
one at $h_i = 0$ that corresponds to the holding cost~\eqref{e-holding-cost}
and one at $h_i = h^{\rm init}_i$ that corresponds to the holding cost~\eqref{e-holding-cost}.
Due to the larger scale of $f_i$ for Johnson \& Johnson,
discontinuities are less severe, relative to the rest of the function,
than they are for the other two assets.
In addition to these discontinuities,
there is also a nonconvex kink present in all three plots around $h^{\rm init}_i$,
due to the tax liability $L_i$ \citep[\S 3]{moehle2020tax}.
Despite these nonconvexities, the separable cost functions 
$f_i$ are generally well approximated
by their convex envelopes $f_i^{**}$,
which helps explain why the upper bound from ADMM and the lower bound
(given in section~\ref{s-results}) are so close.

\clearpage

\printbibliography

\appendix

\section{Portfolio construction parameters}
\label{s-portfolio-params}
We can convert problem~\eqref{e-app-port-orig}
into the separable--affine form~\eqref{e-sap}.
The variable is $x = (h, c, y)\in\reals^{l+k+1}$.
The affine constraint parameters are
\begin{align*}
A = \begin{bmatrix} C^TX & 0 & -I \\ \ones^T & 1 & \phantom{-}0 \end{bmatrix},
\qquad
b = \begin{bmatrix} 0 \\ 1  \end{bmatrix}.
\end{align*}
The part of the separable function corresponding to asset $i$ is
\[
f_i(h_i) =
\alpha_ih_i - \gamma^{\rm risk}D_{ii} h_i^2 + \gamma^{\rm trd} \phi^{\rm trd}_i(h_i - h^{\rm init}_i)
+ \gamma^{\rm hld} \phi^{\rm hld}_i(h_i)
\]
for $i = 1,\dots, l$.
The component of $f$ corresponding to $c$ is
\[
f_{l+1}(c) = \begin{cases}
0 & \text{if $1 - \eta^{\rm ub} \le c \le 1 - \eta^{\rm lb}$} \\
\infty & \text{otherwise.}
\end{cases}
\]
The component of $f$ corresponding to $y$ is
$f_i(y_i) = y_i^2$, for $i = l+2,\dots, l+k+1$.

\section{Convex envelope details}
\label{s-envelope-appendix}
We can compute the convex envelope of a piecewise quadratic function
as discussed in section~\ref{s-cvx-env-pwq}.
To do this, we must compute the parameters
$\alpha$, $\beta$, $x_\varphi$, and $x_\psi$ from \eqref{e-env-two-cvx-funs}.
We do this making use of the following two facts.

The first fact is that candidate values of $\alpha$, $\beta$, $x_\varphi$, and $x_\psi$
are valid if and only if the function $\alpha x + \beta$ is tangent to $\psi^i$ at $x_\psi$
and to $\varphi_i$ at $x_\varphi$.
This holds if and only if
\begin{align}
\label{e-env-param-cond-pwq}
\alpha \in \partial \psi^i(x_\psi), \quad
\psi^i(x_\psi) = \alpha x_\psi + \beta, \quad
\alpha \in \partial \varphi_i(x_\varphi), \quad
\varphi_i(x_\varphi) = \alpha x_\psi + \beta,
\end{align}
where $\partial \psi^i$ and $\partial \varphi_i$ are the subdifferentials
of $\psi^i$ and $\varphi_i$.
These conditions can be easily checked given a set of candidate parameters
$\alpha$, $\beta$, $x_\psi$, and $x_\varphi$.

The second fact is that \eqref{e-env-two-cvx-funs} can be rewritten
\begin{equation}
  \label{e-env-two-cvx-funs-mod}
  \psi^{i+1}(x)
  = \min\{\psi^i, \varphi_i\}^{**}
  = \begin{cases}
    \psi_1^i(x) & \text{if $x \in [\tilde a_1, \tilde b_1]$} \\
    \hfill\vdots\hfill & \hfill\vdots\hfill  \\
    \psi_{j-1}^i(x) & \text{if $x \in [\tilde a_{j-1}, \tilde b_{j-1}]$} \\
    \min\{\psi_j^i, \varphi_i \}^{**} & \text{if $x \in [\tilde a_j, b_i]$,} \\
  \end{cases}
\end{equation}
where $\tilde a_1 \le \tilde b_i \le \dots \le \tilde a_j$ are the endpoints of the pieces of $\psi^i$.
In other words, the function $\psi^{i+1}$
matches $\psi^i$ up to the $j$th piece of $\psi^i$,
and after that, is equal to the convex envelope of
the pointwise minimum of $\psi^i_j$ and $\varphi_i$.
This convex envelope (the last piece in \eqref{e-env-two-cvx-funs-mod})
is easy to compute, as $\min\{\psi_j^i, \varphi_i\}$
is a piecewise quadratic with two pieces.
(We discuss how to do this in section~\ref{impl:env:two-piece}.)

This means that once $j$ is known, we can use \eqref{e-env-two-cvx-funs-mod}
to compute $\psi^{i+1}$.
To find $j$, we simply try all pieces of $\psi^i$,
compute the right-hand side of \eqref{e-env-two-cvx-funs-mod},
and check if the resulting parameters $\alpha$, $\beta$, $x_\psi$ and $x_\varphi$
satisfy the conditions \eqref{e-env-param-cond-pwq}.

\subsection{Piecewise-quadratic functions with two pieces}
\label{impl:env:two-piece}
The function $\min\{\psi_j^i, \varphi_i\}$ in \eqref{e-env-two-cvx-funs-mod}
is piecewise quadratic with two pieces.
Here we discuss how to compute the convex envelope of such functions.

Let $g(x)$ be a PWQ with two pieces: $g_1(x) = p_1x^2 + q_1x + r_1$ on $[a_1, b_1]$, and
$g_2(x) = p_2x^2 + q_2x + r_2$ on $[a_2, b_2]$, with $b_1 \leq a_2$.
In this case, $g^{**}$ has the form
\begin{equation}
  \label{enveq}
  g^{**}(x) = \begin{cases}
      g_1(x) & \text{if $x \in[a_1, x_1]$}\\
      h(x) = \alpha x + \beta & \text{if $x \in [x_1, x_2]$}\\
      g_2(x) & \text{if $x \in [x_2, b_2]$}
  \end{cases}
\end{equation}
for some $\alpha$, $\beta$, $x_1$, and $x_2$.
These parameters are real valued,
but we allow for the case when $x_1 = -\infty$ or $x_2 = \infty$;
in these cases, first or last interval is degenerate, and can be ignored.

Similarly to \eqref{e-env-param-cond-pwq} above,
it is necessary and sufficient for the
parameters $\alpha$, $\beta$, $x_1$, and $x_2$ to satisfy
\begin{align}
\label{e-env-param-cond}
\alpha \in \partial g_1(x_1), \quad
g_1(x_1) = \alpha x_1 + \beta, \quad
\alpha \in \partial g_2(x_2), \quad
g_2(x_2) = \alpha x_2 + \beta.
\end{align}
How these checks are carried out in practice
depends on whether $x_1$ (or $x_2$) are in the interior or boundary 
of the domain of $g_1$ (or $g_2$),
or whether $x_1 = -\infty$ (or $x_2 = \infty$).

\subsubsection{The midpoint-to-midpoint case}
\label{impl:env:case:mm}
We first consider the case when $x_1$ is in the interior of the domain of $g_1$,
\ie, $a_1 < x_1 < b_1$, and $x_2$ is in the interior of the domain of $g$, \ie,
$a_2 < x_2 < b_2$. In this case, we must have
\begin{align}
  g_1(x_1) = h(x_1), \quad g_2(x_2) = h(x_2),
  \quad g_1'(x_1) = h(x_1),
  \quad g_2'(x_2) = h(x_2). \label{e-mm}
\end{align}
By plugging in the values of the functions and their derivatives, we obtain
\begin{align}
\label{e-mm-explicit}
\begin{split}
  p_1x_1^2 + q_1x_1 + r_1 &= \alpha x_1 + \beta \\
  p_2x_2^2 + q_2x_2 + r_2 &= \alpha x_2 + \beta \\
  2p_1x_1 + q_1 &= \alpha \\
  2p_2x_2 + q_1 &= \alpha. \\
\end{split}
\end{align}
These four equations can be reduced to a single quadratic equation with a single unknown.
This quadratic equation has at most two solutions, each corresponding to a set
of candidate values of $\alpha$, $\beta$, $x_1$, and $x_2$. To see if these candidate
values $\alpha$, $\beta$, $x_1$, and $x_2$ parameterize a valid convex envelope
of $g$, we check if our initial assumption, that $a_1 < x_1 < b_1$
and $a_2 < x_2 < b_2$, holds.

\subsubsection{The midpoint-to-endpoint case}
\label{impl:env:case:me}
Now we consider the case in which $x_1$ is in the interior of the domain of $g_1$,
and $x_2 \in \{a_2, b_2\}$. (We note that the case in which $x_1 \in \{a_1, b_1\}$
and $x_2$ in the interior of the domain of $g$ can be handled similarly, and we
do not discuss it further.)

\paragraph{Finite upper bound.}
We start with the case when $x_2 < \infty$.
To do this, we solve a slightly modified set of equations:
\begin{align}
\label{e-me}
  g_1(x_1) = h(x_1), \quad g_1'(x_1) = h'(x_1), \quad g_2(x_2) = h(x_2).
\end{align}
These are the first three equations of \eqref{e-mm-explicit}.
As before, they can be reduced to a single quadratic equation with a single unknown.
This quadratic equation has at most two solutions,
each corresponding to a set of candidate values of $\alpha$, $\beta$, and $x_1$.

To see if these candidate values parameterize a valid convex envelope,
we verify that $x_1$ is in the interior of the domain of $g_1$, \ie, $a_1 < x_1 < b_1$,
and also the second condition of \eqref{e-env-param-cond}.
In the degenerate case in which $g_2$ is defined over a single point,
\ie, $a_2 = b_2$, this second condition always holds;
in the non-degenerate case $a_2 < b_2$,
the condition is equivalent to $g_2'(x_2) = 2p_2x_2 + q_2 \le \alpha$ if $x_2 = a_2$,
and $g_2'(x_2) = 2p_2x_2 + q_2 \ge \alpha$ if $x_2 = b_2$.

\paragraph{Infinite upper bound.}
Next we consider the case of $x_2 = \infty$, which may occur when $x_2 = b_2 = \infty$,
\ie, the domain of $g_2$ is unbounded.
In this case, we modify the last equation in \eqref{e-me}, resulting in the equations
\begin{align}
\label{e-me-mod}
  g_1(x_1) = h(x_1), \quad g_1'(x_1) = h'(x_1), \quad g_2'(\infty) = h'(\infty).
\end{align}
The last equation is equivalent to $p_2 = 0$ and $\alpha = q_2$.
To find the corresponding candidate values of $\beta$, and $x_1$,
we solve $g_1(x_1) = h(x_1)$ and $g_1'(x_1) = h'(x_1)$.
To check the validity of the candidate values of $\alpha$, $\beta$, $x_1$, and $x_2$,
we check that $a_1 < x_1 < b_1$ and $p_2 = 0$.

\subsubsection{The endpoint-to-endpoint case}
\label{impl:env:case:ee}
Finally, we consider the case in which $x_1 \in \{a_1, b_1\}$ and $x_2\in \{a_2, b_2\}$.

\paragraph{Finite upper bound.}
First suppose $a_1$, $b_1$, $a_2$, and $b_2$ are all finite.
Take $h$ to be the line through the points $(x_1, g_1(x_1))$ and $(x_2, g_2(x_2))$,
\ie, the parameters are $\alpha = (g(x_2) - g_1(x_1))/(x_2 - x_1)$
and $\beta = \alpha x_1 - g_1(x_1)$. Note that if $x_1 = x_2$, then $h$ is ill-defined,
and these candidate values of $x_1$ and $x_2$ can be skipped.

To verify that these values of $\alpha$, $\beta$, $x_1$, and $x_2$
parameterize a valid convex envelope of $f$, we check condition~\eqref{e-env-param-cond}.
Recall that the first condition only need hold if $a_1 < b_1$,
and the second if $a_2 < b_2$. For example, if $a_1=b_1<a_2=b_2$,
then any combination of $x_1 \in \{a_1, b_1\}$ and $x_2 \in \{a_2, b_2\}$
are immediately valid.

\paragraph{Infinite upper bound.}
Now we consider the case when $x_2 = b_2 = \infty$.
In this case, we require $h'(\infty) = g_2'(\infty)$,
\ie, $p_2 = 0$ and $\alpha = q_2$. We then have $\beta = g_1(x_1) - q_2 x_1$.
To verify that $\alpha$, $\beta$, $x_1$, and $x_2$ parameterize a valid envelope
of $f$, we need only check the first condition of \eqref{e-env-param-cond}.

\end{document}